\documentclass[11pt] {article}

\usepackage{a4}
\usepackage [francais] {babel}

\begin{document}


\newcommand\tens[1]{\mathop\otimes_{^#1}}

\newcommand\tilda[1]{{\vbox{\hbox{$\sim$}\vskip -6pt plus 0pt minus 0pt\hbox{$#1$}}}}
\newcommand\rond[1]{{\displaystyle\mathop{#1}^o}}


\newcounter {ItemUn}
\newcounter {ItemDeux}
\newcounter {ItemTrois}
\newcounter {ItemNb}
\newcounter {ListeNb}
\newcounter {soussection}
\newcounter {cont}
\newcounter {Chapter}
\newcounter {contref}
\newcounter {notePageNllePage}
\newcounter {notePageNum}
\setcounter {notePageNum} {1}


\def\lastref{}

\newwrite\fileref
\newif\ifrefmodif

\def\initfileref{%
    \input\jobname.ref%
    \refmodiffalse%
    \immediate\openout\fileref=\jobname.ref%
    \global\let\initfileref=\relax}%

\def\defineref#1#2{{\def\next{#1}%
    \expandafter\xdef\csname ref = \meaning\next\endcsname{#2}%
    }}

\def\label#1{{%
    \toks0={#1}\wlog{REF \the\toks0=\lastref}%
    \initfileref
    \immediate\write\fileref{\noexpand\defineref{\the\toks0}{\lastref}}%
    \def\next{#1}%
    \expandafter\ifx\csname ref = \meaning\next\endcsname\lastref
                \else{\global\refmodiftrue
                      \message{{ATTENTION : la reference \next  a change}}
                     }\fi
    \defineref{#1}{\lastref}%
    }}

\def\ref#1{{%
    \initfileref\def\next{#1}%
    \csname ref = \meaning\next\endcsname}}

\outer\def\Ebye{%
    \closeout\fileref%
    \ifrefmodif\erreurmodif\fi%
    }

\def\erreurmodif{%
    \wlog{ATTENTION, certaines r\'{e}f\'{e}rences ont chang\'{e}es.}%
    \wlog{Veuillez recompiler le texte.}
    \message{ATTENTION, certaines references ont changees.}%
    \message{Veuillez recompiler le texte.}}


\newcommand \ligne {\par\small{~}\par}
\def\etape(#1,#2){{
    \par
    \vskip 5pt plus 0pt minus 0pt
    \noindent\textit{\underbar{Step #1}~: #2}\par}}
\def\point(#1,#2){{\par\noindent\textit{\underbar{Point #1}~: #2}\par}}

\newcommand\interconteur{
    \stepcounter{cont}
    \global\edef\lastref{\arabic{section}.\arabic{cont}}
    }
\newcommand\enviraux[2]{%
    \bfseries
    \interconteur
    {\tiny{~}\par}
    \noindent#1 \arabic{section}.\arabic{cont}.~#2 %
    \normalfont}%
\newcommand\envirauxx[2]{%
    \bfseries
    {\tiny{~}\par}
    \noindent#1#2. %
    \normalfont}%

\newenvironment {preuve} [1] [] {\envirauxx{Proof} {#1}} {}
\newenvironment {proposition} [1] [] {\enviraux{Proposition} {#1}} {}
\newenvironment {note} [1] [] {\enviraux{#1} {}} {}
\newenvironment {theoreme} [1] [] {\enviraux{Theorem} {#1}} {}
\newenvironment {definition} [1] [] {\enviraux{Definition} {#1}} {}
\newenvironment {remarque} [1] [] {\enviraux{Remark} {#1}} {}
\newenvironment {exemple} [1] [] {\enviraux{Example} {#1}} {}
\newenvironment {lemme} [1] [] {\enviraux{Lemme} {#1}} {}


\def\chapdebut{
         \par
         \cleardoublepage
         \stepcounter{Chapter}%
         \setcounter{section}{0}%
         \setcounter{soussection}{0}%
         \setcounter{cont}{0}%
         \global\edef\lastref{Chap. \Roman{Chapter}}%
         {\Large{~}\par}{\Large{~}\par}%
         \hrule height 1pt depth 1pt width \hsize
         \vskip 15pt plus 0pt minus 0pt
         }
\def\chapchap{
         \centerline{\large Chapitre \arabic{Chapter}}
         \vskip 10pt plus 0pt minus 0pt
         }

\newcommand\chapter[2][]{
         \chapdebut
         \chapchap
         \centerline{\LARGE\bfseries\itshape
                \vbox{\hbox to \hsize{\hss #1\hss}
                      \vskip -5pt plus 0pt minus 0pt
                      \hbox to \hsize{\hss #2\hss}}}
         \vskip 15pt plus 0pt minus 0pt
         \hrule height 1pt depth 1pt width \hsize
         \vskip 40pt plus 0pt minus 0pt
         }

\newcommand\chaptersn[2][]{
         \chapdebut
         \centerline{\LARGE\bfseries\itshape
                \vbox{\hbox to \hsize{\hss #1\hss}
                      \vskip -5pt plus 0pt minus 0pt
                      \hbox to \hsize{\hss #2\hss}}}
         \vskip 15pt plus 0pt minus 0pt
         \hrule height 1pt depth 1pt width \hsize
         \vskip 40pt plus 0pt minus 0pt
         }

\def\section#1{
         \stepcounter{section}%
         \setcounter{soussection}{0}%
         \setcounter{cont}{0}%
         \global\edef\lastref{\arabic{section}.}%
         \vskip 40pt plus 0pt minus 0pt
         {\noindent%
         \large\bfseries\arabic{section}.~#1}%
         \par
         \vskip 20pt plus 1pt minus 1pt
         }

\def\subsection#1{
         \stepcounter{soussection}
         \setcounter{cont}{0}
         \global\edef\lastref{\Roman{section}.\alph{soussection}.}
         \vskip 15pt plus 1pt minus 1pt
         {\parindent 2cm\large\itshape\underbar{\alph{soussection})~#1}}
         \par
         \vskip 8pt plus 1pt minus 1pt
}


\newfam\msbmfam
\font\tenmsbm=msbm10 \textfont\msbmfam=\tenmsbm
\font\sevenmsbm=msbm7 \scriptfont\msbmfam=\sevenmsbm
\font\fivemsbm=msbm5 \scriptscriptfont\msbmfam=\fivemsbm
\def\msbm{\fam\msbmfam\tenmsbm}
\def \SetR {{\msbm R}}
\def \SetZ {{\msbm Z}}
\def \SetN {{\msbm N}}


\def\textindent#1{\indent\llap{#1\enspace}\ignorespaces}
\newenvironment{ListeNum}
    {\vskip 0,2cm\addtocounter{ItemNb}{1}\setcounter{ListeNb}{1}}
    {\ifcase\value{ItemNb}
        {}
        \or{\setcounter{ItemUn}{0}}
        \or{\setcounter{ItemDeux}{0}}
        \or{\setcounter{ItemTrois}{0}}
    \fi\addtocounter{ItemNb}{-1}\vskip 0,2cm}

\newenvironment{Liste}
    {\addtocounter{ItemNb}{1}\setcounter{ListeNb}{2}}
    {\ifcase\value{ItemNb}
        {}
        \or{\setcounter{ItemUn}{0}}
        \or{\setcounter{ItemDeux}{0}}
        \or{\setcounter{ItemTrois}{0}}
    \fi\addtocounter{ItemNb}{-1}\par}

\def\Item{\par\vskip 4pt plus 0pt minus 0pt%
    \ifcase\value{ItemNb}%
        \or\parindent=1.5cm\hangindent=\parindent\addtocounter{ItemUn}{1}%
            \ifcase\value{ListeNb}%
                    \or\textindent{\arabic{ItemUn}.}
                    \or\textindent{$\bullet$}
                    \fi
        \or\parindent=3cm\hangindent=\parindent\addtocounter{ItemUn}{1}%
            \ifcase\value{ListeNb}%
                    \or\textindent{\arabic{ItemUn}.}
                    \or\textindent{-}
                    \fi
        \or\parindent=4.5cm\hangindent=\parindent\addtocounter{ItemUn}{1}%
            \ifcase\value{ListeNb}%
                    \or\textindent{\arabic{ItemUn}.}
                    \or\textindent{$\ast$}
                    \fi
    \fi}


\def\hhsize{15,5cm} \hsize=\hhsize \hoffset=-0,9cm
\def\hsizesurdeux{3cm}
\renewcommand{\thepage}{\hbox to \hsizesurdeux{\hss- \arabic{page} -}}


\newcommand\Enote [1]{\if%
    {notePageNllePage=\count0}{}\else%
    {\setcounter{notePageNllePage}{0}}\setcounter{notePageNum}{1}\fi%
    \footnote [\value{notePageNum}] {\hsize=\hhsize\noindent#1}}%

\author{E. Reynaud}
\begin{center}
{~}\vskip 1cm
{\LARGE Algebraic fundamental group and simplicial complexes}%
{\vskip 0.5cm \large Eric Reynaud$^*$}\insert\footins{\hskip
-8,7cm\footnotesize
$^*$ E-mail address : reynaud@math.univ-montp2.fr\hss}%
{\vskip 0.25cm D\'{e}partement de Math\'{e}matiques, Universit\'{e} de Montpellier II}%
{\par pl. Eug\`{e}ne Bataillon, 34095 Montpellier Cedex 5\par} \vskip
1cm
\end{center}

\hrule \vskip0,5cm\noindent\textbf{Abstract} \vskip0,3cm

In this paper we prove that the fundamental group of a simplicial
complex is isomorphic to the algebraic fundamental group of its
incidence algebra, and we derive some applications.

\vskip0,5cm\noindent{\small  \textit{AMS classification} : 16E40;
16G20; 06A11; 55Q05}

\vskip0,5cm\hrule\vskip1,3cm

Let $k$ be a field and $A$ be a basic and split finite dimensional
k-algebra, which means that $A/r=k\times k\times\ldots \times k$
where $r$ is the radical of $A$. There exists a unique quiver $Q$
and usually several admissible ideals $I$ of the algebra $kQ$ such
that $A=kQ/I$ (see [6]). In the 1980s, an algebraic fundamental
group has been defined which depends on the presentation of $A$,
that is to say on the choice of the ideal $I$ (see [13]). For
incidence algebras, that is algebras obtained from a simplicial
complex, it has been proved that the presentation does not
influence the fundamental group ([15]). Then it is a natural
question to compare it with the fundamental group of the geometric
realisation. Note also that in [4,8] the analogous question
concerning homology is solved.

Actually, we prove that the fundamental groups considered for a
finite and connected simplicial complex are isomorphic. The
following diagram summarizes the situation:

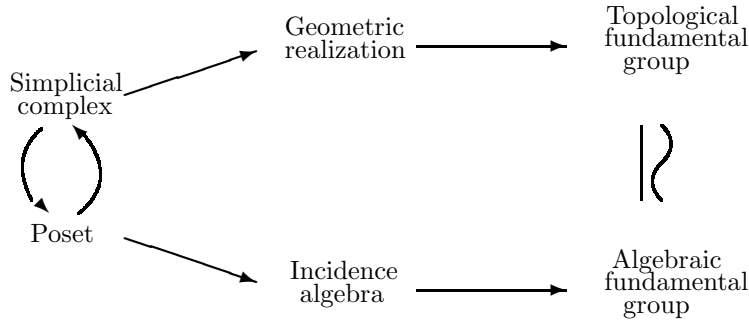
\begin{center}\begin{picture} (300,170)(0,-20)
\unitlength=0.25mm \thicklines%
\put(14,130){\small Simplicial}\put(18,118){\small complex}%
\put(160,160){\small Geometric}\put(160,148){\small realization}%
\put(331,166){\small Topological}%
\put(330,154){\small fundamental}%
\put(340,142){\small group}%
\put(25,50){\small Poset}%
\put(163,30){\small Incidence}%
\put(168,18){\small algebra}%
\put(334,36){\small Algebraic}
\put(330,24){\small fundamental} \put(340,12){\small group}%
\put(75,128){\vector(3,1){70}} \put(75,52){\vector(3,-1){70}}%
\put(230,154){\vector(1,0){80}} \put(230,24){\vector(1,0){80}}%
\bezier{100}(25,72)(13,95)(30,110)\put(25,75){\vector(1,-1){10}}
\bezier{100}(50,65)(72,83)(52,108)\put(57,102){\vector(-1,1){10}}
\put(350,72){\line(0,1){40}} \bezier{100}(360,72)(350,82)(360,92)
\bezier{100}(360,92)(370,102)(360,112)
\end{picture}\end{center}

The isomorphism above enables us to adapt results of algebraic
topology to the purely algebraic setting : for instance the
isomorphism recently proved between $Hom(\Pi_1(Q,I_{_Q}),k^+)$ and
$H\!H^1(kQ/I_{_Q})$ where $Q$ is an ordered quiver, $I_{_Q}$ the
associated parallel path and $k^+$ the additive group of a field
([8], see also [7]) is a consequence of the classic result in
algebraic topology which states that the abelianisation of the
$\Pi_1$ and the first homology group of a simplicial complex are
isomorphic. In another direction, we also derive an algebraic Van
Kampen theorem.

For our purpose we consider in the first part of this paper
simplicial complexes and usual fundamental groups ; as well as
posets, incidence algebras and algebraic fundamental groups. We
also explain the relations between posets and simplicial
complexes. In the last section, we provide the isomorphism and
give applications.

This work is part of my thesis in Montpellier and I thank C.
Cibils for his help, his apposite advice and his patience. When
finishing this paper, I learned that J.C. Bustamente had
independently considered a similar context.

All simplicial complexes will be finite, connected and not empty.


\section{Fundamental group and incidence algebras}

%
%
We consider the classical definition of a simplicial complex (see
for instance [11]). A simplicial complex is the union of elements
$\{a_i\}_{i\in I}$ called the vertices, and finite sets of
vertices called the simplexes, such that if  $S$ is a simplex all
subsets of $S$ are also simplexes. Each set containing only one
element is a simplex.
A simplicial complex is said to be finite if $I$ is finite, and is
connected if for each couple of vertices $(s,t)$, there exist
vertices $s_0,\ldots,s_n$ such that $s_0=s$, $s_n=t$ and
$\{s_{i-1},s_i\}$ is a simplex for each $i$ in $\{1,\ldots,n\}$.
%

The geometric realization of a simplicial complex $C$ with
vertices $\{a_i\}_{i\in I}$ is as follows. Let $\{A_i\}_{i\in I}$
be points of $\SetR^n$, such that if
$\{a_{\alpha_1},\ldots,a_{\alpha_p}\}$ is a simplex of $C$, then
the points $A_{\alpha_1},\ldots,A_{\alpha_p}$ are linearly
independent. The set of points whose barycentric coordinates are
strictly positive is called a face. Note that we prefer to define
a face as the points with strictly positive barycentric
coordinates instead of positive barycentric coordinates: this way,
the geometric realization becomes the disjoint union of its faces
and there exists only one face containing it. Moreover, if $S$ and
$S'$ are two simplexes with empty intersection, the intersection
of the corresponding faces is empty. A geometric realization of
$C$ denoted by $|C|$, is the union of the faces associated to
simplexes of a simplicial complex. A closed face is the closure of
a face.

Note that since we assume that the simplicial complexes are
finite, their geometric realization exist. Indeed let $C$ be a
complex with set of vertices $a_1,\ldots,a_n$.
%
%
%
%
%
%
%
The geometric realization $|C|$ of a simplicial complex $C$ is a
subset of $\SetR^n$ and inherits the topology of $\SetR^n$. Behind
the usual definition of the fundamental group $\Pi_1(|C|)$
obtained through homotopy classes of closed paths, we recall the
construction of the edge-paths group of $|C|$. This provides
another description of the fundamental group $\Pi_1(|C|)$ which is
useful for our purpose (see [11] for example).

%
%
%

An \textit{edge-path} of $C$ is a finite sequence of vertices
$a_{i_r}\ldots a_{i_1}$ such that for all $j$ in
$\{1,\ldots,r-1\}$, the set $\{a_{i_j},a_{i_{j+1}}\}$ is a
simplex. If $w=a_{i_r}\ldots a_{i_1}$ is an edge-path, let
$w^{-1}$ denote the edge-path $a_{i_1}\ldots a_{i_r}$. An
edge-path is said to be \textit{closed} (or \textit{an edge-loop})
if the first and the last vertices are the same. Let
$w=a_{i_r}\ldots a_{i_1}$ and $w'=a'_{i'_{r'}}\ldots a'_{i'_1}$ be
two edge-paths; if $a_{i'_{r'}}=a_{i_1}$, the product $w.w'$ is
defined and is equal to $a_{i_r}\ldots a_{i_1}a'_{i'_{r'}}\ldots
a'_{i'_1}$

This is an \textit{allowable operation} on edge-paths : if three
consecutive vertices of the edge-path are in the same simplex, the
middle vertex can be removed. Conversely, we can add a vertex
between two others, if these three vertices are in a same simplex
of $C$. Moreover, it is possible to change $a_{i_0}a_{i_0}$ by
$a_{i_0}$ and conversely. This generates an equivalence relation
on the set of edge-paths. Let $\overline{w}$ denote the
equivalence class of the edge-path $w$. As two equivalent
edge-paths have the same extremities, the product defined before,
when it exists, is defined also on the equivalence classes, as
well as on the set of edge-loops starting at a fixed point.

\begin{proposition} [ ({[11] 6.3.1 and 6.3.2.})]
    Let $C$ be a simplicial complex and $x_0$ a vertex of $C$.
    The set of equivalence classes of edge-loops starting at one point
    $x_0$ is a group for the product defined before.
    Since $C$ is connected, this group does not depend on $x_0$
    and is denoted by $\Pi_1(C)$.
    The fundamental groups
    $\Pi_1(C)$ and $\Pi_1(|C|)$ are isomorphic.
\end{proposition}\vskip0,2cm

Hereafter, $\Pi_1(C)$ will be either the approximation of the
fundamental group $\Pi_1(C)$ or the fundamental group $\Pi_1(|C|)$
itself.





Given a \textit{poset} (i.e. a partially ordered set), there is an
associated \textit{ordered quiver}, that is to say a finite
oriented graph without loops and such that if there exists an
arrow from $a$ to $b$, there does not exist another path from $a$
to $b$. To each element of the poset corresponds a vertex of the
graph. Moreover, let $S_1$ and $S_2$ be vertices in the graph ;
there exists an arrow from $S_1$ to $S_2$ if and only if the
element associated to $S_1$ in the poset is smaller than the
element associated to $S_2$ and if there does not exist an element
of the poset strictly between these two elements. The graph
obtained is an ordered quiver.

For example, the graph that corresponds to the poset
${a,b,c,a',b',c',d}$ with $a<b'<d$, $c<b'<d$, $a<c'<d$, $b<c'<d$,
$b<a'<d$ and $c<a'<d$ is:

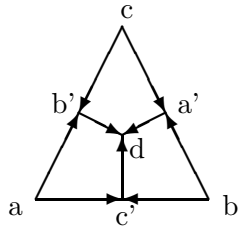
\begin{center}\begin{picture} (100,90) (0,10)
\unitlength=0.23mm%
\thicklines%
\put(20,20){\vector(1,0){50}}
        \put(120,20){\vector(-1,0){50}} \put(20,20){\vector(1,2){25}}
        \put(120,20){\vector(-1,2){25}} \put(70,120){\vector(-1,-2){25}}
        \put(70,120){\vector(1,-2){25}} \put(45,70){\vector(2,-1){25}}
        \put(95,70){\vector(-2,-1){25}} \put(70,20){\vector(0,1){37}}

        \put(4,11){a} \put(66,5){c'} \put(128,11){b}
        \put(29,70){b'} \put(69,125){c}
        \put(102,70){a'} \put(74,44){d}
\end{picture}\end{center}
\vskip 20pt plus 0pt minus 0pt

Conversely, by this operation all ordered graphs arise from a
poset. There is a bijection between the set of ordered graphs and
the set of posets. Moreover, a poset is said to be connected if
its ordered quiver is connected. All the posets considered will be
connected, finite and non empty.



Let now $Q$ be a quiver, and $k$ be a field. We denote $kQ$ the
$k$-vector space with basis the paths of $Q$ (the paths of length
0 being the vertices), with the multiplication given by the
composition of two paths if possible and 0 otherwise.
Two paths of $Q$ are \textit{parallel} if they have the same
beginning and the same end. The k-space generated by the set of
differences of two parallel paths is a two-sided ideal of $kQ$,
denoted $I_Q$ and called \textit{parallel ideal}.
The quotient algebra $kQ/I_Q$ is the \textit{incidence algebra} of
$Q$.



The general definition of the fundamental group depends on a
couple $(Q,I)$ where $I$ is an admissible ideal of $kQ$, it can be
found in [8, 15, 7] and we recall it below. We notice that for an
algebra the presentation as a quiver with relations is not unique
in general, that is to say that different ideals $I$ may exist
such that $A\cong kQ/I$. We do not have in general a unique
fundamental group associated to an algebra. Nevertheless, in the
case of an incidence algebra it has been proved that the
fundamental group does not depend on the presentation of the
algebra (see [15] for example).

A relation $\sum_{i=1}^n \lambda_i\omega_i$ is \textit{minimal} if
the sum is in $I$ and if for all non empty proper subset $J$ of
$\{1,\ldots,n\}$ the sum $\sum_{i\in J} \lambda_i\omega_i$ is not
in $I$. We note that if the relation is minimal then
$\{\omega_1,\ldots,\omega_n\}$ have the same source and the same
terminus.

If $\alpha$ is an arrow from $x$ to $y$, let $\alpha^{-1}$ denote
its formal inverse which goes from $y$ to $x$. A walk from $x$ to
$y$ is a formal product
$\alpha_1^{\pm1}\alpha_2^{\pm1}\ldots\alpha_n^{\pm1}$ which begins
in $x$ and ends in $y$. The trivial walk $x$, which begins in $x$
and no longer moves is denoted $e_x$. A \textit{closed walk} (or a
\textit{loop}) is a walk having the same extremities. A walk is,
in fact, a path in the non oriented graph associated to $Q$; in
other words, a walk can follow the arrows in any direction.

We consider the smaller equivalence relation $\sim$ on the walks
of $Q$ containing the following items :
    \begin{ListeNum}
        \Item if $\alpha$ is an arrow from $x$ to $y$ then
        $\:\:\alpha\alpha^{-1} \sim e_y\:\:$ and $\:\:\alpha^{-1}\alpha \sim
        e_x$,
        \Item if $\sum_{i=1}^n \lambda_i\omega_i$ is a minimal relation
        then $\omega_1\sim\ldots\sim\omega_n$,
        \Item if $\;\;\alpha\sim\beta\;\;$ then, for all
        $(\omega,\omega')$, we have
        $\omega\:\alpha\:\omega'\sim\omega\:\beta\:\omega'$
    \end{ListeNum}
Let $x_0$ be a vertex of $Q$. The set of equivalence classes of
loops starting at $x_0$ does not depend on $x_0$ since $Q$ is
connected. We denote this set by $\Pi_1(Q,I)$. If the quiver $Q$
comes from a poset $P$, the fundamental group associated to
$(Q,I_Q)$ is denoted $\Pi_1(P)$.

We remark that if $I$ is the parallel ideal, the second item means
parallel paths are equivalent. For example if $Q$ is defined by
the following quiver, the fundamental group $\Pi_1(Q,I_{_Q})$ is
isomorphic to~$\SetZ$.

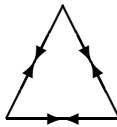
\begin{center}\begin{picture} (60,40) (0,5)
\unitlength=0.15mm \thicklines \put(20,20){\vector(1,0){50}}
\put(120,20){\vector(-1,0){50}} \put(20,20){\vector(1,2){25}}
\put(120,20){\vector(-1,2){25}} \put(70,120){\vector(-1,-2){25}}
\put(70,120){\vector(1,-2){25}}
\end{picture}\end{center}

We note that the hypothesis $I$ admissible is not fully used
neither to define the fundamental group nor to prove that it is a
group. So, in this paper, we will consider the fundamental group
of a couple $(Q,I)$ where $I$ satisfies $F^n\subset I\subset F$
for an integer $n$. This will be used to adapt Van Kampen's
theorem to purely algebraic fundamental groups.





Let $C$ be a simplicial complex. The set of non empty simplexes of
$C$ ordered by inclusion is a poset which we will denote $Pos(C)$.
This $Pos$ defines an application from simplicial complexes to
posets which is injective, but not surjective, since there is no
simplicial complex which gives the poset {$a<b$}.
%
%
%
Due to the construction of the quiver from a simplicial complex,
an arrows can only go from a vertex corresponding to a p-face to a
vertex corresponding to a q-face with $p>q$, then the path algebra
of this quiver is then of finite dimension.




We provide now the construction of a simplicial complex from a
poset. These procedures are of course not inverse one of each
other, their composition is the barycentric decomposition.

To each poset $P$ we associate a simplicial complex $Sim(P)$,
where a n-simplex is a subset of $P$ containing $n+1$ elements and
totally ordered. The application $Sim$ is surjective but not
injective, for instance the simplicial complexes which are
associated to $a<b<c$, $a<b<d$ and to $a<c<b$, $a<d<b$ are the
same.




Let $C$ be a simplicial complex, $|Sim(Pos(C))|$ is the geometric
realization of the barycentric decomposition of $C$. Then for
example let $C$ be the set of the non empty parts of $\{a,b,c\}$ ;
its geometric realization being the triangle drawn on the next
figure. Then, $Pos(C)$ contains all the elements of $C$ that is to
say $T=\{a,b,c\}$ $A_1=\{b,c\}$, $A_2=\{a,c\}$, $A_3=\{a,b\}$,
$S_1=\{a\}$, $S_2=\{b\}$ and $S_3=\{c\}$ and its order is defined
by $S_i\leq A_j\leq T$ for all $i,j\in\{1,2,3\}$ and $i\neq j$.
The associated quiver is drawn on the next figure. Then
$Sim(Pos(C))$ is the complex containing the total ordered subsets
of $P$ that is to say : for all $i,j\in\{1,2,3\}$, the 0-simplexes
are $\{S_i\}$, $\{A_i\}$, $\{T\}$, the 1-simplexes are
$\{S_i,A_j\}$ $i\neq j$, $\{S_i,T\}$, $\{A_i,T\}$, the 2-simplexes
are $\{S_i,A_j,T\}$ $i\neq j$. By identifying poset and associated
ordered quiver, simplicial complexes and their geometric
realization, the situation can be summarized by the following
diagram :

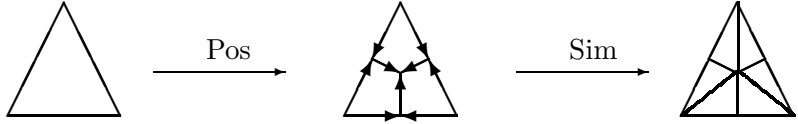
\begin{center}\begin{picture} (100,100) (0,-30)
        \unitlength=0.15mm\thicklines \put(20,20){\line(1,0){50}}
        \put(120,20){\line(-1,0){50}} \put(20,20){\line(1,2){25}}
        \put(120,20){\line(-1,2){25}} \put(70,120){\line(-1,-2){25}}
        \put(70,120){\line(1,-2){25}}
\end{picture}
\begin{picture}(10,0)(40,-55)\put(0,0){\vector(1,0){50}}
        \put(20,5){Pos}\end{picture}
\begin{picture} (100,100) (-10,-30)
        \unitlength=0.15mm\thicklines \put(20,20){\vector(1,0){50}}
        \put(120,20){\vector(-1,0){50}} \put(20,20){\vector(1,2){25}}
        \put(120,20){\vector(-1,2){25}} \put(70,120){\vector(-1,-2){25}}
        \put(70,120){\vector(1,-2){25}} \put(45,70){\vector(2,-1){25}}
        \put(95,70){\vector(-2,-1){25}} \put(70,20){\vector(0,1){37}}
\end{picture}
\begin{picture}(10,0)(20,-55)\put(0,0){\vector(1,0){50}}
        \put(20,5){Sim}\end{picture}
\begin{picture} (100,100) (-20,-30)
        \thicklines
        \unitlength=0.15mm \put(20,20){\line(1,0){50}}
        \put(120,20){\line(-1,0){50}} \put(20,20){\line(1,2){25}}
        \put(120,20){\line(-1,2){25}} \put(70,120){\line(-1,-2){25}}
        \put(70,120){\line(1,-2){25}} \put(45,70){\line(2,-1){25}}
        \put(95,70){\line(-2,-1){25}} \put(70,20){\line(0,1){37}}
        \bezier{50}(20,20)(45,40)(70,60)\bezier{50}(120,20)(95,40)(70,60)
        \bezier{50}(70,120)(70,90)(70,60)
\end{picture}\end{center}



%



\vskip -1,5cm\section{Equivalence between algebraic and
topological approaches.}

The aim of this section is to prove that the fundamental group
$\Pi_1(|C|)$ defined on the geometric realization of a finite
simplicial complex $C$ is isomorphic to the fundamental group
$\Pi_1(Pos(C))$ of the incidence algebra of the poset deduced from
the complex.

We prove first that for any poset $P$ we have
$\Pi_1(P)\simeq\Pi_1(Sim(P))$. We will use the approximation of
the topological fundamental group considered before %
in order to provide an isomorphism
between $\Pi_1(Sim(P))$ and $\Pi_1(P)$.


\begin{theoreme}
Let $P$ be a poset, the groups $\Pi_1(P)$ and $\Pi_1(Sim(P))$ are
isomorphic. 
\end{theoreme}

\begin{preuve}
    Due to proposition 1, it is sufficient to prove that $\Pi_1(P)$
    is isomorphic to the edge-paths group of $|Sim(C)|$

    For each vertex $s$ of the poset $P$, the set $\{s\}$ is totally
    ordered and it corresponds to a 0-simplex of the simplicial complex
    $Sim(P)$. We denote by $s'$ this 0-simplex of $Sim(P)$.

    Let $\phi$ be the map from
    the set of the associated quiver walks to the set of edge-paths of
    $Sim(P)$ be defined by
    $$\phi(\alpha_n^{\epsilon_n}\ldots\alpha_1^{\epsilon_1})=s'_{n+1}\ldots s'_1$$
    where $s_i$ and $s_{i+1}$ are the poset elements which are the origin and the end of the walk
    $\alpha_i^{\epsilon_i}$.


    \vskip0,3cm
    The map $\phi$ is well defined.
    Indeed, $s'_{n+1}\ldots s'_1$ is an edge-path because
    for all $i$ in $\{1,\ldots,n\}$,
    the set $\{s_i,s_{i+1}\}$ is totally ordered since $\alpha_i$ is an arrow
    of the extremities $s_i,s_{i+1}$ and
    therefore
    $\{s'_i,s'_{i+1}\}$ is a simplex of $Sim(P)$.

    \vskip0,3cm
    We assert now that the images by $\phi$ of
    equivalent walks are equivalent.
    We have to prove this fact on the generators of the
    equivalence relation. Let $f$ be an arrow from $s_1$ to
    $s_2$, then $\phi(f.f^{-1})=s'_1.s'_2.s'_1$ which is equivalent to
    $s'_1=\phi(e_{v_1})$. Let $c_1$ and $c_2$ be two parallel paths crossing
    respectively the vertex $s_1,t_1,t_2,\ldots,t_n,s_2$ and
    $s_1,u_1,u_2,\ldots,u_n,s_2$. Then the sets
    $\{s_1,t_1,t_2,\ldots,t_n,s_2\}$ and
    $\{s_1,u_1,u_2,\ldots,u_n,s_2\}$ are totally ordered and the
    edge paths $\phi(c_1)=s'_1.t'_1.t'_2.\ldots.t'_n,s'_2$ and $\phi(c_2)=s'_1.t'_1.t'_2.\ldots.t'_n,s'_2$
    are both equivalent to $s'_1.s'_2$.

    The third relation is immediate because the equivalence
    relation on the edge paths set is compatible with the product of the group.

    Let $\phi_*:\Pi_1(P)\rightarrow \Pi_1(Sim(P))$ denote the
    application induced by $\phi$. Note that the image of a loop is also a loop.

    First of all $\phi_*$ is a morphism. Indeed :

    \noindent If
    $\left\{\begin{array} {l}
    \phi(p)=x_0 v_n\ldots v_1 x_0 \\
    \phi(q)=x_0 w_n\ldots w_1 x_0 \\
    \end{array}\right.$
    then
    $\left\{\begin{array} {l}
    \phi(p).\phi(q)=x_0 v_n\ldots v_1x_0x_0 w_n\ldots w_1 x_0 \\
    \phi(p.q)=x_0 v_n\ldots v_1 x_0 w_n\ldots w_1 x_0\\
    \end{array}\right.$

    \noindent These two edge-paths are equivalent so $\phi_*(p)\phi_*(q)=\phi_*(p.q)$.

    \vskip0,3cm To prove that $\phi_*$ is bijective, we are going
    to construct its inverse $\psi_*$.

    Let $s'_{n+1}\ldots s'_1$ be an edge-path. We fix $i$ in $\{1,\ldots,n\}$.
    The set $\{s'_i,s'_{i+1}\}$ being a simplex, the set $\{s_i,s_{i+1}\}$
    is totally ordered. Then, there exists a maximal (for the inclusion) totally ordered set
    containing it and having $s_i$ and $s_{i+1}$ as extremities.
    The choice of this is not important because all paths
    corresponding to these sets have the same origin and the same
    end and therefore are parallel.
    This maximal set
    corresponds to a path or to an inverse path $w_i^{\epsilon_i}$
    of
    the associated quiver with origin $s_i$ and end $s_{i+1}$.
     So we can define a morphism
    $\psi$ from edge paths group to $\Pi_1(C)$ by
    $\psi(s'_{n+1}\ldots s'_1)=w_n^{\epsilon_n}\ldots
    w_1^{\epsilon_1}$ if $n\geq 1$ and $\psi(s'_1)=s_1$.

    We will prove now that this application
    is constant on the equivalence class. In deed $\psi(s's')$ is a
    path from $s$ to $s$, so $\psi(s's')=s=\psi(s')$. Moreover,
    let take $s',t',u'$ such that $\{s',t',u'\}$ is a simplex, so $\psi(s't'u')=
    \psi(s't').\psi(t'u')$ and $\psi(s'u')$ are paths from $u$ to
    $s$. Therefore they are parallel.

    Finally, we will verify that
    $\phi_*\:o\:\psi_*=\psi_*\:o\:\phi_*=Id$. Let $f$ be an arrow
    of the ordered quiver associated to $P$ from $s_1$ to $s_2$
    then $\phi_*(f)=s_1's_2'$ and $\psi_*\:o\:\phi_*(f)$ is a path
    from $s_1$ to $s_2$. Since the quiver is ordered, it does not
    exist any other path than $f$, so $\psi_*\:o\:\phi_*(f)=f$ and
    $\psi_*\:o\:\phi_*=Id$.
    For the other equality, let's consider an edge path
    $s_1's_2'$. Then $\psi_*(s_1's_2')$ is a path from $s_1$ to
    $s_2$ and then $\phi_*\:o\:\psi_*(s_1's_2')$ is an edge path
    beginning with $s_1'$ and ending with $s_2'$ such that all
    vertices of this edge-path are in a same simplex. So it is
    equivalent to $s_1's_2'$.
\end{preuve}


\begin{theoreme}
    Let  $C$ be a simplicial complex. The fundamental groups $\Pi_1(|C|)$
    and $\Pi_1(Pos(C))$ are isomorphic.
\end{theoreme}

\begin{preuve}
    The fundamental groups associated to the simplicial complex
    and to its barycentric decomposition are
    isomorphic. Then, $$\Pi_1(|C|)\simeq\Pi_1(Sim(Pos(C)).$$
    Moreover, the previous theorem shows that
    the group $\Pi_1(Pos(C))$ and $\Pi_1(Sim(Pos(C)))$ are isomorphic.
\end{preuve}

\section{Applications}
We first show that in a particular case of incidence algebra the
result obtained by I. Assem and J.A. De La Pe\~{n}a ([1], p.200) is a
consequence of a classic fact in algebraic topology ; so this
proof has the advantage to link this result to algebraic topology.

\begin{theoreme}
Let $P$ be a poset and $A=kQ/I$ its incidence algebra over a field
$k$. Then $$Hom(\Pi_1(P),k^+)\simeq H\!H^1(A)$$ where $k^+$ is the
additive group of the field $k$ and $H\!H^1(A)$ is the first
Hochschild cohomology group of $A$.
\end{theoreme}







\begin{preuve}
    Let $C$ be a simplicial complex and $C_\bullet(\Lambda)$ be the
    chain complex over a ring $\Lambda$ associated to the simplicial
    homology of $C$ ; the latter homology will be denoted by
    $H_\bullet(C,\Lambda)$. Moreover, the
    cohomology of the complex
    $Hom_\Lambda(C_\bullet(\Lambda),A)$, where $A$ is a
    $\Lambda$-module, is denoted by $H^\bullet(C,\Lambda,A)$.

    For any simplicial complex, the abelianisation $\Pi_1^{ab}$ of $\Pi_1$ is
    isomorphic to the first homology group (see for example [11],
    6.4.7). If we denote by $C$ the simplicial complex associated to $Q$,
    we have
    $$(\Pi_1(C))^{ab} \simeq H_1(C,\SetZ)$$
    and
    $$Hom((\Pi_1(C))^{ab},k^+)\simeq Hom(H_1(C,\SetZ),k^+).$$
    Note that since $k^+$ is abelian we have :
    $$Hom(\Pi_1(C),k^+)\simeq Hom(H_1(C,\SetZ),k^+).$$
    Making use of theorem %
    we have :
    $$Hom(\Pi_1(Q,I_{_Q}),k^+)\simeq Hom(H_1(C,\SetZ),k^+)\:\:\:\:\:\:(*)$$

    Adjunction, for two rings $R$ and $S$ and bimodules
    $A_R$, $_R\!B_S$, $C_S$, provides an isomorphism, (see [17], p.37, for
    example) :
    $$Hom_S(A\tens {R}B,C)\simeq Hom_R(A,Hom_S(B,C))$$
    In our case, we consider $B=C=S=k$, $R=\SetZ$ and
    $A=H_1(C,\SetZ)$. Identifying commutative groups and $\SetZ$-modules, this
    isomorphism becomes
    $$Hom_k(H_1(C,\SetZ)\tens {\SetZ}k,k)\simeq Hom(H_1(C,\SetZ),k^+).$$
    Then $(*)$ becomes
    $$Hom(\Pi_1(Q,I_{_Q}),k^+)\simeq Hom_k(H_1(C,\SetZ)\tens {\SetZ}k,k)\:\:\:\:\:\:(*)$$

    Moreover, the universal coefficients theorem in homology and
    in cohomology (see for example [10], p.176-179) are as
    follows.
    Let $\Lambda$ be a principal ring and
    $C$ be a flat chain complex over $\Lambda$, and let $A$ be a
    $\Lambda$-module, then

    $$\left\{
        \begin{array}{l}
            0\rightarrow H_n(C)\tens{\Lambda} A
                \rightarrow H_n(C\tens{\Lambda} A)
                \rightarrow Tor_1^\Lambda(H_{n-1}(C),A)
                \rightarrow 0\\
            0\rightarrow Ext_\Lambda^1(H_{n-1}(C),A)
                \rightarrow H^n(Hom_\Lambda(C,A))
                \rightarrow Hom_\Lambda(H_n(C),A)
                \rightarrow 0\\
        \end{array}\right.$$ are exact.

    As $C_\bullet(\SetZ)$ is free and therefore flat, the
    hypotheses are verified with
    $C=C_\bullet(\SetZ)$, $\Lambda=\SetZ$, $A=k$ and $n=1$ in the first theorem
    and $C=C_\bullet(k)$, $\Lambda=k$, $A=k$ et $n=1$ in the second.

    Moreover,
    $Ext^1_k(H_{n-1}(C),A)=0$ because $k$ is a field, and
    $Tor_1^\Lambda(H_0(C),k)=0$ because $H_0(\SetZ)$ is free (see for example [10], p.63, cor.2.4.7).

    Then we have
    $$\left\{
    \begin{array}{l}
      H_1(C,\SetZ)\tens{\SetZ} k \simeq H_1(C\tens{\SetZ} k,\SetZ) \\
      H^1(C,k,k))\simeq Hom_k(H_1(C,k),k) \\
    \end{array}
    \right.$$

    Since the complexes
    $C_\bullet(\SetZ)\otimes_\SetZ k$ and $C_\bullet(k)$ are isomorphic, the equation $(*)$ becomes
    $$Hom(\Pi_1(Q,I),k^+)
            \simeq Hom_k(H_1(C\tens{\SetZ} k,\SetZ),k)
            \simeq Hom_k(H_1(C,k),k)
            \simeq H^1(C,k,k)$$
    We conclude using the Gerstenhaber-Schack theorem
    ([4] 1.4, or [8]) which shows that the cohomologies $H^i(C,k,k)$
    and ${H\!H}^i(kQ/I)$ are isomorphic.
\end{preuve}




\vskip0,8cm We introduce now the notion of completed quiver which
is interesting in view of the following results, and also in order
to adapt Van Kampen's theorem to an algebraic setting. A quiver is
said to be a completed quiver if
\begin{ListeNum}
    \Item it does not contain cycles.
    \Item each path of length at least two has a parallel arrow,
    \Item there are no couples of parallel arrows.
\end{ListeNum}

\begin{theoreme}
    Let $Q$ be a quiver without cycles and without parallel arrows.
    \begin{ListeNum}
        \Item There exists an ordered quiver $Q^o$ obtained
        by considering the set of paths of length at least 2 and
        deleting all arrows parallel to these paths.
        \Item There exists an completed quiver $Q^c$ obtained
        from $Q$ by considering the set of paths of length at
        least 2 and adding a parallel arrow for each of these paths,
        unless there already is one, either added or in $Q$.
        \Item The three fundamental groups of the quivers $Q$, $Q^o$ and $Q^c$ with their
        own parallel ideals which are not necessary admissible, are isomorphic:
            $$\Pi_1(Q,I_{_Q})\simeq\Pi_1(Q^o,I_{_{Q^o}})\simeq\Pi_1(Q^c,I_{_{Q^c}})$$
    \end{ListeNum}
\end{theoreme}

\begin{preuve}

First, we give more details on the construction of $Q^c$ and
$Q^o$. The vertices of $Q^c$ and $Q^o$ are the same as those of
$Q$. Moreover, there exists an arrow from $a$ to $b$, $a\neq b$,
in $Q^c$ if and only if there exists a path (possibly an arrow)
from $a$ to $b$ in $Q$, and there exists an arrow in $Q^o$ if and
only if there exists an arrow in $Q$ from $a$ to $b$ which is not
parallel to another path in $Q$.%
\vskip0,3cm%
We prove that for any path in $Q^c$ (resp. in $Q$) there exists a
parallel path constructed with arrows that come from $Q$ (resp.
$Q^o$). Let $q=\alpha_1\ldots\alpha_n$ be a path of $Q^c$. For any
$\alpha_i$ that does not come from $Q$, there exists by
construction a path $\omega_i$ in $Q^c$, already present in $Q$,
parallel to $\alpha_i$. If $\alpha_i$ is an arrow which comes from
$Q$, let us set $\omega_i=\alpha_i$. Then the path $q$ is parallel
to the path $q'=\omega_1\ldots\omega_n$ which is issued from $Q$.

The same process is applicable to the quivers $Q$ and $Q^o$.
\vskip0,3cm%
The fundamental groups $\Pi_1(Q,I_{_Q})$, $\Pi_1(Q^o,I_{_{Q^o}})$
and $\Pi_1(Q^c,I_{_{Q^c}})$ are isomorphic. Since any path in
$Q^c$ is parallel to a path that comes from $Q$, every walk of
$\Pi_1(Q^c,I_{_{Q^c}})$ is equivalent to a walk that contains only
arrows issuing from $Q$ ; then, the equivalent classes will not be
changed by adding the parallel arrows of $Q^c$.

The proof of the isomorphism between $\Pi_1(Q^o,I_{_{Q^o}})$ and
$\Pi_1(Q,I_{_{Q}})$ is the same. %
\vskip0,3cm%
$Q^o$ is ordered. The quiver $Q^o$ does not contain cycles, and
there is no arrow parallel to a path by construction.
\vskip0,3cm%

$Q^c$ is completed. The quiver $Q^o$ does not contain cycles
either, and if there exists a path $p$ from $a$ to $b$ and no
arrow from $a$ to $b$, then there is a path in $Q$, parallel to
$p$. This is in contradiction with the construction of $Q^o$.
\end{preuve}

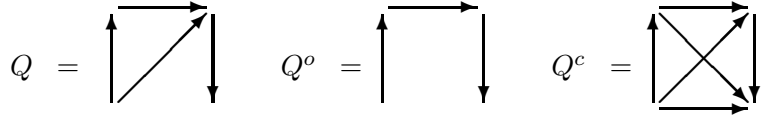
\begin{exemple}
\begin{center}\begin{picture}(330,60)(0,0)
    \thicklines
    \unitlength=0.45mm
    \put(10,20){$Q\;\;=$}
    \put(40,12){\vector(0,1){26}}
    \put(42,40){\vector(1,0){26}}
    \put(70,38){\vector(0,-1){26}}
    \put(42,12){\vector(1,1){26}}
    \put(90,20){$Q^o\;\;=$}
    \put(120,12){\vector(0,1){26}}
    \put(122,40){\vector(1,0){26}}
    \put(150,38){\vector(0,-1){26}}
    \put(170,20){$Q^c\;\;=$}
    \put(200,12){\vector(0,1){26}}
    \put(202,40){\vector(1,0){26}}
    \put(230,38){\vector(0,-1){26}}
    \put(202,12){\vector(1,1){26}}
    \put(202,38){\vector(1,-1){26}}
    \put(202,10){\vector(1,0){26}}
\end{picture}\end{center}
\end{exemple}

\begin{theoreme} [ (Van Kampen)]
    Let $Q$ be a connected ordered quiver and let $Q_1$, $Q_2$ be two
    connected subquivers such that $Q_1^c\cup Q_2^c=Q^c$ and such that $Q_0$,
    the ordered quiver associated to the
    completed quiver $Q_1^c\cap Q_2^c$, is also connected. The inclusions $i_1$
    (resp. $i_2$) from $Q^c_0$ to $Q^c_1$ (resp. $Q^c_2$)
    induce $i_{1*}$ (resp $i_{2*}$) from $\Pi_1(Q_0,I_{_{Q_0}})$ to $\Pi_1(Q_1,I_{_{Q_1}})$
    (resp. $\Pi_1(Q_2,I_{_{Q_2}})$).

    Therefore $\Pi_1(Q,I_{_{Q}})$ is the free product of $\Pi_1(Q_1,I_{_{Q_1}})$ and
    $\Pi_1(Q_2,I_{_{Q_2}})$  by adding the relations
    $i_{1*}(\alpha)=i_{2*}(\alpha)$ for all $\alpha$ in $\Pi(Q_0,I_{_{Q_0}})$
\end{theoreme}

\begin{preuve}
We use Van Kampen's theorem for a simplicial complex (see [11],
p243, for example):

\itshape{\leftskip=1.5cm\rightskip=1.5cm $K$ is a connected
complex and $K_0$, $K_1$, $K_2$ are connected subcomplexes, such
that $K_1\cup K_2=K$ and $K_1\cap K_2=K_0$; $a_0$ is a vertex of
$K_0$; the injection maps $j_r:K_0\rightarrow K_r$, $r=1,2$,
induce $$j_{r*}:\Pi(K_0,a_0)\mapsto\Pi(K_r,a_0)$$}
{\leftskip=1.5cm\rightskip=1.5cm Then $\Pi_1(K)$ is the free
product of $\Pi_1(K_1)$ and $\Pi_1(K_2)$ by means of the mapping
$j_{1*}(\alpha)\mapsto j_{2*}(\alpha)$ for all $\alpha$ in
$\Pi(Q_0,I_{_{Q_0}})$.\par}\normalfont

\vskip 0,3cm \noindent To $Q$,$Q_0$,$Q_1$,$Q_2$, there correspond
simplicial complexes that we denote $K$,$K_0$,$K_1$,$K_2$. By
hypothesis, these simplicial complexes are connected. Moreover, as
there is a bijection between the paths of $Q^c_i$ of length $n$
and the $n$-dimensional simplexes of $K_i$, we have $K=K_1\cup
K_2$ and $K_0=K_1\cap K_2$.

Let us denote by $\phi_*$ the isomorphism defined in theorem
between $\Pi_1(Q,I_{_{Q}})$ and $\Pi_1(K)$.
Then it is easy to see that the following diagrams are commutative
(n=1,2):

\begin{center}\begin{picture} (150,100)(0,0)
\put(20,20){$\Pi_1(K_0)$} \put(0,80){$\Pi_1(Q_0,I_{_{Q_0}})$}
\put(110,20){$\Pi_1(K_n)$} \put(110,80){$\Pi_1(Q_n,I_{_{Q_n}})$}
\put(55,23){\vector(1,0){50}} \put(55,83){\vector(1,0){50}}
\put(40,70){\vector(0,-1){35}} \put(130,70){\vector(0,-1){35}}
\put(20,50){$\phi_*$} \put(140,50){$\phi_*$}
\put(75,10){$j_{k*}$}\put(75,90){$i_{k*}$} 
\bezier{100}(80,63)(95,63)(95,53)
\bezier{100}(95,53)(95,43)(80,43) \put(80,43){\vector(-1,0){1}}
\end{picture}\end{center}

Therefore $\Pi_1(Q,I_{_{Q}})$ is the free product of
$\Pi_1(Q_1,I_{_{Q_1}})$ and $\Pi_1(Q_2,I_{_{Q_2}})$  by adding the
relations $i_{1*}(\alpha)=i_{2*}(\alpha)$ for all $\alpha$ in
$\Pi(Q_0,I_{_{Q_0}})$
\end{preuve}

\begin{note} [Corollary]
If $\Pi_1(Q_0,I_{_{Q_0}})=0$, then $\Pi_1(Q,I_{_{Q}})$ is the free
product of $\Pi_1(Q_1,I_{_{Q_1}})$ and $\Pi_1(Q_2,I_{_{Q_2}})$.
\end{note}

\begin{exemple}
    We begin with an easy example to show how the theorem works
    \noindent

    \begin{center}$Q=$
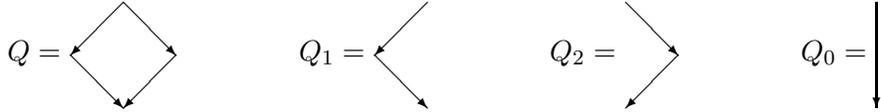
    \begin{picture} (40,40)(0,18)
    \put(20,40){\vector(-1,-1){20}}
    \put(20,40){\vector(1,-1){20}}
    \put(0,20){\vector(1,-1){20}}
    \put(40,20){\vector(-1,-1){20}}
    \end{picture}
    \hskip1.5cm
    $Q_1=$
    \begin{picture} (20,40)(0,18)
    \put(20,40){\vector(-1,-1){20}}
    \put(0,20){\vector(1,-1){20}}
    \end{picture}
    \hskip1.5cm
    $Q_2=$
    \begin{picture} (20,40)(0,18)
    \put(0,40){\vector(1,-1){20}}
    \put(20,20){\vector(-1,-1){20}}
    \end{picture}
    \hskip1.5cm
    $Q_0=$
    \begin{picture} (10,40)(0,18)
    \put(0,40){\vector(0,-1){40}}
    \end{picture}
    \end{center}

    \vskip 1,5cm
    \noindent As the fundamental group of a tree is zero, we deduce from the corollary that
    $\Pi_1(Q,I_{_{Q}})=0$.
\end{exemple}

\begin{exemple}
    Now, let's consider the quiver $Q$ drawn on the next figure. The
    dotted line indicates the limits of quivers $Q_1$ and $Q_2$.
    Then the quivers $Q$, $Q_1$, $Q_2$ and $Q_0$ are :

    
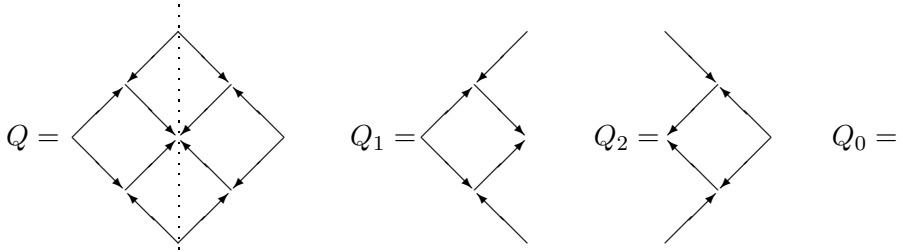
\begin{center}\begin{picture}(100,100)(0,0)
  {
    \put(-5,37){$Q=$}
    \bezier{25}(60,90)(60,40)(60,-10)
    \put(20,40){\vector(1,1){19}}
    \put(40,60){\vector(1,-1){19}}
    \put(20,40){\vector(1,-1){19}}
    \put(40,20){\vector(1,1){19}}
    \put(60,80){\vector(1,-1){19}}
    \put(60,80){\vector(-1,-1){19}}
    \put(60,0){\vector(1,1){19}}
    \put(60,0){\vector(-1,1){19}}
    \put(100,40){\vector(-1,1){19}}
    \put(100,40){\vector(-1,-1){19}}
    \put(80,60){\vector(-1,-1){19}}
    \put(80,20){\vector(-1,1){19}}
  }
    \end{picture}
    \hskip1cm\begin{picture}(60,100)(0,0)
  {
    \put(-7,37){$Q_1=$}
    \put(20,40){\vector(1,1){19}}
    \put(40,60){\vector(1,-1){19}}
    \put(20,40){\vector(1,-1){19}}
    \put(40,20){\vector(1,1){19}}
    \put(60,80){\vector(-1,-1){19}}
    \put(60,0){\vector(-1,1){19}}
  }
    \end{picture}
    \hskip1cm
    \begin{picture}(60,100)(0,0)
    {
    \put(-7,37){$Q_2=$}
    \put(20,80){\vector(1,-1){19}}
    \put(20,0){\vector(1,1){19}}
    \put(60,40){\vector(-1,1){19}}
    \put(60,40){\vector(-1,-1){19}}
    \put(40,60){\vector(-1,-1){19}}
    \put(40,20){\vector(-1,1){19}}
    }
    \end{picture}
    \hskip1cm
    \begin{picture}(30,100)(0,0)
    {
    \put(-9,37){$Q_0=$}
    \put(20,80){\vector(0,-1){39}}
    \put(20,0){\vector(0,1){39}}
    }
    \end{picture}\end{center}

    Since the quiver $Q_1$ and $Q_2$ are the same, we only had to
    calculate the fundamental group of $Q_1$. Let's use again the
    Van Kampen theorem, and decompose the quiver $Q_1$ in $Q_{11}$
    and $Q_{12}$. The intersection quiver will be denoted by
    $Q_{10}$ :

    
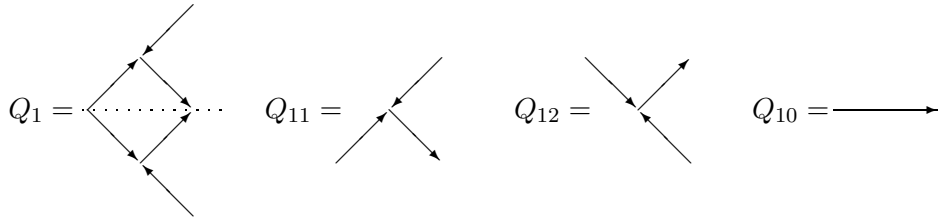
\begin{center}
    \begin{picture}(60,100)(0,0)
        {\put(-10,37){$Q_1=$}
         \bezier{12}(18,40)(44,40)(70,40)
         \put(20,40){\vector(1,1){19}}\put(40,60){\vector(1,-1){19}}%
         \put(20,40){\vector(1,-1){19}}\put(40,20){\vector(1,1){19}}%
         \put(60,80){\vector(-1,-1){19}}\put(60,0){\vector(-1,1){19}}}%
    \end{picture}\hskip1.2cm
    \begin{picture}(60,100)(0,0)
        {\put(-7,37){$Q_{11}=$}
         \put(20,20){\vector(1,1){19}}\put(40,40){\vector(1,-1){19}}%
         \put(60,60){\vector(-1,-1){19}}}%
    \end{picture}\hskip1.2cm
    \begin{picture}(60,100)(0,0)
        {\put(-7,37){$Q_{12}=$}
         \put(20,60){\vector(1,-1){19}}\put(40,40){\vector(1,1){19}}%
         \put(60,20){\vector(-1,1){19}}}%
    \end{picture}\hskip1.2cm
    \begin{picture}(60,100)(0,0)
        {\put(-11,37){$Q_{10}=$}
         \put(20,40){\vector(1,0){39}}}%
    \end{picture}\hskip1.2cm
    \end{center}
    The fundamental groups of quivers $Q_{10}$, $Q_{11}$ and $Q_{12}$ is
    0. Then it is the same for the quiver $Q_1$ and therefore for $Q$ itself.
\end{exemple}

\begin{exemple}
    The last example shows a situation where $Q_0$ is not simply
    connected :
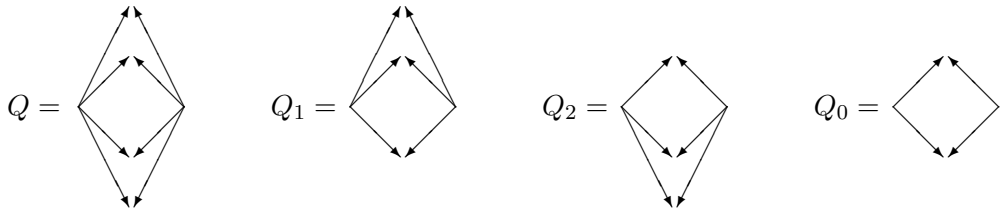
\begin{center}
\begin{picture}(60,100)(0,0)
  {
    \put(-7,57){$Q=$}
    \put(20,60){\vector(1,1){19}}\put(20,60){\vector(1,2){19}}%
    \put(20,60){\vector(1,-1){19}}\put(20,60){\vector(1,-2){19}}%
    \put(60,60){\vector(-1,1){19}}\put(60,60){\vector(-1,2){19}}%
    \put(60,60){\vector(-1,-1){19}}\put(60,60){\vector(-1,-2){19}}%
  }
\end{picture}\hskip1.5cm
\begin{picture}(60,100)(0,0)
  {
    \put(-10,57){$Q_1=$}
    \put(20,60){\vector(1,1){19}}\put(20,60){\vector(1,2){19}}%
    \put(20,60){\vector(1,-1){19}}%
    \put(60,60){\vector(-1,1){19}}\put(60,60){\vector(-1,2){19}}%
    \put(60,60){\vector(-1,-1){19}}%
  }
\end{picture}\hskip1.5cm
\begin{picture}(60,100)(0,0)
  {
    \put(-10,57){$Q_2=$}
    \put(20,60){\vector(1,1){19}}%
    \put(20,60){\vector(1,-1){19}}\put(20,60){\vector(1,-2){19}}%
    \put(60,60){\vector(-1,1){19}}%
    \put(60,60){\vector(-1,-1){19}}\put(60,60){\vector(-1,-2){19}}%
  }
\end{picture}\hskip1.5cm
\begin{picture}(60,100)(0,0)
  {
    \put(-10,57){$Q_0=$}
    \put(20,60){\vector(1,1){19}}%
    \put(20,60){\vector(1,-1){19}}%
    \put(60,60){\vector(-1,1){19}}%
    \put(60,60){\vector(-1,-1){19}}%
  }
\end{picture}\end{center}
Once again, the quivers $Q_1$ and $Q_2$ are the same and we use
again the Van Kampen theorem to calculate it :

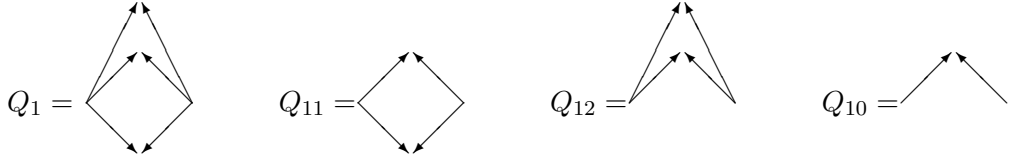
\begin{center}
\begin{picture}(60,100)(0,0)
  {
    \put(-10,57){$Q_1=$}
    \put(20,60){\vector(1,1){19}}\put(20,60){\vector(1,2){19}}%
    \put(20,60){\vector(1,-1){19}}%
    \put(60,60){\vector(-1,1){19}}\put(60,60){\vector(-1,2){19}}%
    \put(60,60){\vector(-1,-1){19}}%
  }
\end{picture}\hskip1.5cm
\begin{picture}(60,100)(0,0)
  {
    \put(-10,57){$Q_{11}=$}
    \put(20,60){\vector(1,1){19}}%
    \put(20,60){\vector(1,-1){19}}%
    \put(60,60){\vector(-1,1){19}}%
    \put(60,60){\vector(-1,-1){19}}%
  }
\end{picture}\hskip1.5cm
\begin{picture}(60,100)(0,0)
  {
    \put(-10,57){$Q_{12}=$}
    \put(20,60){\vector(1,1){19}}\put(20,60){\vector(1,2){19}}%
    \put(60,60){\vector(-1,1){19}}\put(60,60){\vector(-1,2){19}}%
  }
\end{picture}\hskip1.5cm
\begin{picture}(60,100)(0,0)
  {
    \put(-10,57){$Q_{10}=$}
    \put(20,60){\vector(1,1){19}}
    \put(60,60){\vector(-1,1){19}}
  }
\end{picture}\hskip1.5cm
\end{center}
Since the fundamental group $Q_{10}$ is 0, $Q_{1}$ is the free
product of $Q_{11}$ and $Q_{12}$ which are isomorphic to $\SetZ$.
Then, if $Q_{11}$ and $Q_{12}$ are generated by $a$ and $b$, $Q_1$
is the group generated by $\{a,b\}$. In the same way, $Q_2$ is
generated by two elements $c$ and $d$. If we add the relation
$i_{1*}(Q_0)-i_{2*}(Q_0)=0$, we obtain that $d=b$ and therefore
the fundamental group of $Q$ is the group generated by
$\{a,b,c\}$.
\end{exemple}

{\vskip 1cm\noindent\large\textbf{References}\vskip 0,5cm}



\hangindent=\parindent\textindent{[1]}  I. \textsc{Assem} and J.A.
\textsc{De La Pe\~{n}a}, \textit{The fundamental groups of a}
\textit{triangular algebra}, Comm. Algebra, 24(1), pp.187-208
(1996).

\hangindent=\parindent\textindent{[2]}  M. \textsc{Berger},
\textit{G\'{e}om\'{e}trie}, Nathan (1990).

\hangindent=\parindent\textindent{[3]}  K. \textsc{Bongartz} and
P. \textsc{Gabriel}, \textit{Covering Spaces in
representation-Theory}, Invent. Math. 65, pp.331-378 (1982).

\hangindent=\parindent\textindent{[4]}  C. \textsc{Cibils},
\textit{Cohomology of incidence algebras and simplicial complex},
J. Pure Appl. Algebra 56 pp.221-232 (1989).

\hangindent=\parindent\textindent{[5]}  C. \textsc{Cibils},
\textit{Complexes simpliciaux et carquois}, %
C.R. Acad. Sci. Paris t.307, Serie I, pp.929-934 (1988).

\hangindent=\parindent\textindent{[6]}  P. \textsc{Gabriel},
\textit{Indecomposable representation II}, Symposia Mathematica II
(Instituto Nazionale di alta Matematica), Roma, pp.81-104 (1973).

\hangindent=\parindent\textindent{[7]}  M. A. \textsc{Gatica} and
M. J. \textsc{Redondo}, \textit{Hochschild cohomology and
fundamental groups of incidence algebras}, to appear in Comm.
Algebra.

\hangindent=\parindent\textindent{[8]}  M. \textsc{Gerstenhaber}
and S.P. \textsc{Schack}, \textit{Simplicial cohomology is
Hochschild cohomology}, J. Pure Appl. Algebra 30 pp.143-156
(1983).

\hangindent=\parindent\textindent{[9]}  E.L. \textsc{Green},
\textit{Graphs withs relations, coverings and group-graded
algebras}, Trans. Amer. Math. Soc. 279 (1983), 297-310.

\hangindent=\parindent\textindent{[10]}  P.J. \textsc{Hilton} and
U. \textsc{Stammbach}, \textit{A course in Homological Algebra},
Springer (1996).

\hangindent=\parindent\textindent{[11]}  P.J. \textsc{Hilton} and
S. \textsc{Wylie}, \textit{An introduction to Algebraic Topology},
Cambridge University press (1967).

\hangindent=\parindent\textindent{[12]}  W.S. \textsc{Massey},
\textit{Algebraic Topology: an Introduction}, Springer-Verlag,
New-York Heidelberg Berlin (1989).

\hangindent=\parindent\textindent{[13]}  R.
\textsc{Martinez-Villa} and J.A. \textsc{De La Pe\~{n}a}, \textit{The
universal cover of a quiver with relations}, J. Pure Appl. Algebra
30, pp.277-292 (1983).

\hangindent=\parindent\textindent{[14]}  J.A. \textsc{De La Pe\~{n}a}
and M. \textsc{Saorin}, \textit{The first Hochschild cohomology
group of an algebra}, preprint.

\hangindent=\parindent\textindent{[15]}  J.A. \textsc{De La Pe\~{n}a},
\textit{On the abelian Galois covering of an algebra}, J. Algebra
102(1) pp.129-134 (1986).

\hangindent=\parindent\textindent{[16]}  M.J. \textsc{Redondo},
\textit{Cohomolog\'{\i}a de Hochschild de \'{a}lgebras de Artin}, preprint.

\hangindent=\parindent\textindent{[17]}  J. \textsc{Rotman},
\textit{An introduction to Homological Algebra}, Academic press,
inc. (1979).

\Ebye
\end{document}